\documentclass[11pt,reqno]{amsart}
\usepackage{amssymb,latexsym}
\usepackage{graphicx}
\usepackage{fancyhdr}

\newcommand{\N}{{\mathbb N}}
\newcommand{\la}{\lambda}

\theoremstyle{plain} \numberwithin{equation}{section}
\newtheorem{thm}{Theorem}[section]

\newtheorem{Proposition}[thm]{Proposition}
\renewenvironment{proof}[1][\proofname]{{\it #1.}}{\hfill$\square$}

\begin{document}

\setcounter{page}{1}

\title[A NOTE ON THE MODES OF THE POISSON DISTRIBUTION OF ORDER $k$]
{A note on the modes of the Poisson distribution of order $k$}
\author{Andreas~N.~Philippou}
\address{Department of Mathematics\\
               University of Patras\\
               Patras 26500\\
               Greece and \vspace{-2.6mm}}
               \address{Technological Educational Institute of Lamia, Lamia, Greece}
               \email{anphilip@math.upatras.gr}

\begin{abstract}
It is shown that the Poisson distribution of order $k$ ($\geq1$) with
parameter $\la$ ($>0$)  has a unique mode $m_{k,\la}=0$ if
$0<\la<2/(k(k+1))$. In addition, $m_{2,\la}=0$ if
$0<\la\leq-1+\sqrt{3}$ and $m_{2,\la}=2$ if $-1+\sqrt{3}\leq\la<1$.
\end{abstract}

\maketitle\vspace{-6mm}

\section{Introduction}
The Poisson distribution of order $k\,(\geq1)$ with parameter
$\la\,(>0)$ has probability mass function (pmf)
\begin{equation}\label{Pok}
f_{k}(x;\lambda)= e^{-k\lambda} \sum
\frac{\lambda^{x_{1}+\cdots+x_{k}}}{x_{1}!\cdots x_{k}!}, \ \
\mbox{for} \ \ x=0, 1, 2, \ldots,
\end{equation}
where the summation is taken over all $k$-tuples
of non-negative integers $x_{1}, x_{2}$, $\cdots$, $x_{k}$ such that
$x_{1}+2x_{2}+\cdots+kx_{k}=x$ [1--9]. It was obtained in \cite{PGP}
as a limit of the negative binomial distribution of order $k$, and
was named so because it reduces to the Poisson distribution with pmf
$f_1(x;\lambda)= e^{-\la}\la^x/x!$ for $k=1$.

Denote by $m_{k,\la}$ the mode(s) of $f_k(x;\la)$, i.e. the value(s)
of $x$ for which  $f_k(x;\la)$ attains its maximum. It is well known
that $f_1(x;\lambda)$ has a unique mode
$m_{1,\la}=\lfloor\la\rfloor$ if $\la\not\in\N$, and two modes
$m_{1,\la}= \la$ and $\la - 1$ if $\la\in\N$. It was established
recently \cite{GPS} that
\begin{equation}\label{Th1}
\lfloor \mu_{k,\la} \rfloor - {k(k+1)}/{2} +1- \delta_{k,1} \leq
m_{k,\lambda}\leq \lfloor \mu_{k,\la} \rfloor, \ \mbox{for} \
\la>0,k\geq1,
\end{equation}
and
\begin{equation}\label{Th2}
m_{k,\lambda}= \mu_{k,\lambda} -\lfloor k/2 \rfloor, \ \mbox{for} \
\lambda\in\mathbb{N}, 2\leq k \leq5,
\end{equation}

\noindent where $\mu_{k,\lambda}=\la k(k+1)/2$ is the mean of the Poisson
distribution of order $k$ \cite{P1983}, $\delta_{k,1}$ is the Kronecker delta,
and $\lfloor u \rfloor$ denotes the greatest integer not exceeding $u \in \mathbb{R}$. 

We presently show, as a consequence of (1.2), that $m_{k,\lambda}=0$
if $0<\la<2/(k(k+1))$ and $k\geq1$ (see Proposition \ref{p:21}). We
also show that $m_{2,\lambda}=0$ if $0<\la\leq-1+\sqrt{3}$ and
$m_{2,\la}=2$ if $-1+\sqrt{3}\leq\la<1$ (see Proposition
\ref{p:22}).

\section{Main results}

In this section we state and prove the following:

\begin{Proposition}\label{p:21}
For any integer $k\geq1$ and $0<\la<2/(k(k+1))$, the Poisson
distribution of order $k$ has a unique mode $m_{k,\lambda}=0$.
\end{Proposition}

\begin{Proposition}\label{p:22}
The Poisson distribution of order $2$ has a unique mode
$m_{2,\lambda}=0$ if $0<\la<-1+\sqrt{3}$. It has two modes
$m_{2,\lambda}=0$ and $2$ if $\la=-1+\sqrt{3}$, and it has a unique
mode $m_{2,\lambda}=2$ if $-1+\sqrt{3}<\la<1$.
\end{Proposition}

\noindent\begin{proof}[Proof of Proposition 2.1] We have
\begin{eqnarray*}\label{Th1}\begin{array}{ll}\vspace{1mm}
0 & \hspace{-2mm}\leq m_{k,\lambda} \; \mbox{for $k \ge 1$, $\la>0$, by the definition of} \
m_{k,\lambda},\\ \vspace{1mm}
  & \hspace{-2mm}\leq \lfloor \mu_{k,\la} \rfloor = \lfloor \la k(k+1)/2 \rfloor,\,
  \mbox{by} (1.2), \\
  & \hspace{-2mm}=0,\; \mbox{since} \ 0<\la
  k(k+1)/2<1 \ \mbox{by the assumption}.
\end{array}\end{eqnarray*}
For $k=1$, the condition $0< \la < 2/(k(k+1))$ is obviously
necessary and sufficient for Proposition \ref{p:21} to hold true.
For $k=2$, however, it is not necessary because of Proposition
\ref{p:22}.
\end{proof} \vspace{2mm}

\noindent\begin{proof}[Proof of Proposition 2.2] The definition of
$m_{2,\lambda}$ and (1.2) imply
\begin{equation*}
0\leq m_{2,\lambda} \leq \lfloor \mu_{2,\lambda} \rfloor = \lfloor
3\la \rfloor \le 2 \ \mbox{for} \ 0<\la<1.
\end{equation*}
Next, using (\ref{Pok}) or recurrence (2.3) of \cite{GPS}, for
$\la>0$ we get
\begin{equation}\label{S21}
f_2(0;\la)=e^{-2\la}, \ f_2(1;\la)= \la e^{-2\la}, \ \mbox{and} \
f_2(2;\la)= (\la+\frac{\la^2}{2}) e^{-2\la}.
\end{equation}
It follows that
\begin{equation*}
f_2(0;\la) > f_2(1;\la) \ \mbox{for} \ 0<\la<1 , \ \mbox{and} \
f_2(1;\la) < f_2(2;\la) \ \mbox{for} \ \la>0.
\end{equation*}
Therefore, in order to obtain $m_{2,\lambda}$ for $0<\la<1$, it
suffices to compare $f_2(0;\la)$ and $f_2(2;\la)$ for $0<\la<1$. By
means of (\ref{S21}),
\begin{equation*}
f_2(2;\la) \leq f_2(0;\la) \ \mbox{if and only if} \
\la+\frac{\la^2}{2}\leq1  \ \mbox{if and only if} \ 0< \la \leq
-1+\sqrt{3},
\end{equation*}
and
\begin{equation*}
f_2(2;\la) \geq f_2(0;\la) \ \mbox{if and only if} \
\la+\frac{\la^2}{2}\geq1  \ \mbox{if and only if} \ -1+\sqrt{3}\leq
\la < 1,
\end{equation*}
which complete the proof of  the proposition.
\end{proof}

\medskip

\noindent MSC2010: 60E05, 11B37, 39B05

\end{document}